\documentclass[11pt]{amsart}
\usepackage{palatino}
\usepackage[all]{xy,xypic}

\usepackage{amsfonts,amsmath,oldgerm,amssymb,amscd}

\newcommand{\ra}{\rightarrow}		
\newcommand{\lra}{\longrightarrow}
\newcommand{\by}[1]{\stackrel{#1}{\ra}}

\newcommand{\surj}{\ra\!\!\!\ra}	

\newcommand{\ol}{\overline}		

\newcommand{\iso}{\by \sim}

\newtheorem{theorem}{Theorem}[section]
\newtheorem{proposition}[theorem]{Proposition}
\newtheorem{lemma}[theorem]{Lemma}
\newtheorem{definition}[theorem]{Definition}
\newtheorem{corollary}[theorem]{Corollary}

	\newcommand{{\br}}{\mbox{$\mathbb R$}}

	\newcommand{\bz}{\mbox{$\mathbb Z$}}

\newcommand{\mm}{\mbox{$\mathfrak m$}}

\newcommand{\Spec}{\text{Spec}}

\newcommand{\dd}{\text{dim}}

\newcommand{\sur}{\twoheadrightarrow}
\newcommand{\bp}{\begin{proposition}}
\newcommand{\ep}{\end{proposition}}
\newcommand{\bl}{\begin{lemma}}
\newcommand{\el}{\end{lemma}}
\newcommand{\bt}{\begin{theorem}}
\newcommand{\et}{\end{theorem}}
\newcommand{\bc}{\begin{corollary}}
\newcommand{\ec}{\end{corollary}}
\newcommand{\bd}{\begin{definition}}
\newcommand{\ed}{\end{definition}}

\def\rmk{\refstepcounter{theorem}\paragraph{{\bf Remark} \thetheorem}}
\def\proof{\paragraph{Proof}}
\def\example{\refstepcounter{theorem}\paragraph{{\bf Example} \thetheorem}}
\def\notation{\paragraph{\bf Notation}}

\def\definition{\refstepcounter{theorem}\paragraph{{\bf Definition} \thetheorem}}
\title{Euler cycles and Mennicke symbols}
\author{Mrinal Kanti Das}
\address{Stat-Math Unit, Indian Statistical Institute, 203 B. T. Road, Kolkata 700108 India}
\email{mrinal@isical.ac.in}
\author{Soumi Tikader}
\address{Stat-Math Unit, Indian Statistical Institute, 203 B. T. Road, Kolkata 700108 India}
\email{tsoumi\_r@isical.ac.in}
\author{Md. Ali Zinna}
\address{School of Mathematical Sciences, NISER  Bhubaneswar (HBNI),
 Odisha 752050, India.}
\email{zinna2012@gmail.com}

\thanks{We sincerely thank S. M. Bhatwadekar for some very useful discussions. The third named author acknowledges the Department of Science and Technology for their INSPIRE research grant.}

\date{\today}


\begin{document}
\maketitle
\section{Introduction}
Let $R$ be a commutative Noetherian ring of (Krull) dimension $d\geq 2$. The group $E_{d+1}(R)$ (the subgroup of $SL_{d+1}(R)$ generated by the elementary matrices) acts on $Um_{d+1}(R)$, the  set of unimodular rows of length $d+1$ over $R$. When $d=2$, Vaserstein \cite[Section 5]{sv} showed that the orbit space 
$Um_3(R)/E_3(R)$ carries the structure of an abelian group. Later, van der Kallen \cite{vk1}
extended this result
to show that   $Um_{d+1}(R)/E_{d+1}(R)$ has an abelian group structure for all $d\geq 2$. 
This group structure is closely related with the higher Mennicke symbols of Suslin (see
\cite{vk1} for an elaboration).

The group $Um_{d+1}(R)/E_{d+1}(R)$ is intimately related to the $d$\textsuperscript{\,th} Euler class group $E^d(R)$ studied by Bhatwadekar-Sridharan 
(see \cite{brs3,dz,vk3,vk4} for details on the  connection between these two groups). 
The idea of the Euler class group was  envisioned by Nori in order to detect the obstruction for a projective $R$-module of rank $d$ to split off a free summand of rank one. Although this \emph{``splitting problem"} was settled by Bhatwadekar-Sridharan quite sometime back in \cite{brs1,brs3}, surprisingly, the Euler class group has not yet lost its relevance. Very recently, in \cite{dtz2}, the current authors have succeeded in computing the structure of
$Um_{d+1}(R)/E_{d+1}(R)$ for smooth affine ${\mathbb R}$-algebras by comparing this group with the Euler class group, and appealing to the structure theorems for $E^d(R)$ available in \cite{brs2} for such rings. 
To facilitate such a comparison, a set-theoretic map $\delta_R:E^d(R)\lra  Um_{d+1}(R)/E_{d+1}(R)$ 
was defined in \cite{dtz2}, based on the formalism developed in \cite{dtz1}, when $R$ is a smooth affine domain of dimension $d$ over an infinite perfect field $k$ of characteristic unequal to $2$.
If $k={\mathbb R}$, it was proved in \cite{dtz2} that $\delta_R$ is a morphism of groups but at that time it was not clear whether $\delta_R$ is a morphism in general. In this article we prove that 
$\delta_R$ is 
indeed a morphism of groups. We believe this morphism will enable further understanding of these two groups better, as it did in \cite{dtz2} when $k={\mathbb R}$. We must remark that in this article we (re)define 
$\delta_R$ in a much simpler manner than \cite{dtz2} (see Section \ref{maps} for details, and in particular, Remark \ref{original}).

In Sections \ref{ECG} and \ref{homotopy} we recall the definitions of the objects involved in this paper. In Section \ref{maps} we define the map  $\delta_R:E^d(R)\lra  Um_{d+1}(R)/E_{d+1}(R)$. In Section \ref{menn} we treat the special case when the group law in
$Um_{d+1}(R)/E_{d+1}(R)$ is Mennicke-like (as this case is simpler and the treatment is entirely different) and prove that $\delta_R$ is a morphism. In Section \ref{general} we treat the general 
case.

\section{Generalities I: The Euler class group}\label{ECG}

\notation
We shall write an ideal generated by $f_1,\cdots,f_{d}$ as $\langle f_1,\cdots,f_{d}\rangle$.

\smallskip

Let $R$ be a smooth affine domain of dimension $d\geq 2$ over an
infinite perfect field $k$. 
Let $B$ be the set of pairs
$(m,\omega_m)$ where $m$ is a maximal ideal of $R$ and $\omega_m
:(R/m)^d\surj m/m^2$. Let $G$ be the free abelian group generated by
$B$.  Let $J=m_{1}\cap \cdots \cap m_r$, where $m_i$ are distinct maximal
ideals of $R$. Any $\omega_{J}:(R/J)^d\surj J/J^2$ induces surjections
$\omega_i :(R/m_i)^d\surj m_i/m_i^2$ for each $i$. We associate
$(J,\omega_J):= \sum_{1}^{r}(m_i,\omega_i)\in G$.  Now,
 Let $S$ be
the set of elements $(J,\omega_J)$ of $G$ for
which $\omega_J$ has a lift to a surjection $\theta: R^d\surj J$ and $H$ be the 
subgroup of $G$ generated by $S$ .  The
Euler class group $E^d(R)$ is defined as $E^d(R):=G/H$.

\smallskip

\rmk
The above definition appears to be  slightly different from the one given in \cite{brs1}. 
However, note that if $(J,\omega_J)\in S$  and if $\ol\sigma\in E_d(R/J)$, then
the element $(J,\omega_J\ol\sigma)$ is also in  $S$. For details, see  \cite[Proposition 2.2]{dz}. 

\smallskip

\bt\label{zero}\cite[4.11]{brs1}
Let $R$ be a smooth affine domain of dimension $d\geq 2$ over an
infinite perfect field $k$. Let $J\subset R$ be a reduced ideal of height $d$ and
$\omega_{J}:(R/J)^d\surj J/J^2$ be a surjection. Then, the following are equivalent:
\begin{enumerate}
\item
The image of $(J,\omega_J)=0$ in $E^d(R)$
\item
$\omega_J$ can be lifted to a surjection $\theta: R^d\sur J$.
\end{enumerate}
\et

\rmk
We shall refer to the elements of the Euler class group as \emph{Euler cycles}.  
An arbitrary element of $E^d(R)$ can be represented by a single Euler cycle 
$(J,\omega_J)$, where $J$ is a reduced ideal of height $d$ and $\omega_J:(R/J)^d\sur J/J^2$ is
a surjection (see \cite[Remark 4.14]{brs1}).

\smallskip

The following notation  will be used in the rest of this article.

\begin{notation}
Let $\dd(R)=d$.
Let $(J,\omega_J)\in E^d(R)$ and $u\in R$ be a unit modulo $J$. Let $\sigma$ be any diagonal matrix
in $GL_d(R/J)$  with determinant $\ol{u}$ (bar means modulo $J$). 
We shall  denote 
the composite surjection 
$$(R/J)^d\stackrel{\sigma}\iso (R/J)^d\stackrel{\omega_J}\sur J/J^2$$
by  $\ol{u}\omega_J$. It is easy to check that the element $(J,\ol{u}\omega_J)\in E^d(R)$ is
independent of $\sigma$ (the key fact used here is that $SL_{d}(R/J)=E_d(R/J)$ as $\dd(R/J)=0$).
\end{notation}

\section{Generalities II: Homotopy orbits}\label{homotopy}
In this article, by  \emph{`homotopy'} we shall  mean \emph{`naive homotopy'}, as defined below.

\bd
Let $F$ be a functor originating from the category of rings to the category of sets. 
For a given ring $R$, two elements
$F(u_0),F(u_1)\in F(R)$ are said to be homotopic if there is an element $F(u(T))\in F(R[T])$ such that  
$F(u(0))=F(u_0)$ and $F(u(1))=F(u_1)$.
\ed

\bd
Let $F$ be a functor from the category of rings to the category of sets. 
Let $R$ be a ring. Consider the equivalence relation on $F(R)$ \emph{generated by} homotopies
(the relation is easily seen to be reflexive and symmetric but is not transitive in general).
The set of equivalence classes will be denoted by $\pi_0(F(R))$. 
\ed

\medskip

\example
Let $R$ be a ring. Two  matrices $\sigma,\tau\in GL_n(R)$ are \emph{homotopic} if there is a matrix
$\theta(T)\in GL_n(R[T])$ such that $\theta(0)=\sigma$ and $\theta(1)=\tau$. 
Of particular interest are the matrices in $GL_n(R)$ which are \emph{homotopic to identity}.

\bd
Recall that $E_n(R)$ is the subgroup of $GL_n(R)$ generated by all
 elementary matrices $E_{ij}(\lambda_{ij})$ (whose diagonal entries are all $1$, $i\not = j$,
and $ij$-th entry is $\lambda_{ij}\in R$).
\ed

\rmk\label{elem}
Any  $\theta\in E_n(R)$ is homotopic to identity. To see this, let 
$\theta=\prod E_{ij}(\lambda_{ij})$.
Define $\Theta(T):= \prod E_{ij}(T\lambda_{ij})$. Then, clearly $\Theta(T)\in E_n(R[T])$
and we observe that $\Theta(1)=\theta$, $\Theta(0)=I_n$.

In this context, we record below a remarkable
result of Vorst.

\bt\label{vorst}\cite[Theorem 3.3]{v}
Let $R$ be a regular ring which is essentially of finite type over a field $k$. Let $n\geq 3$ and 
$\theta(T)\in GL_n(R[T])$ be such that  $\theta(0)=I_n$ ($\theta$ is thus a homotopy between $I_n$ and $\theta(1)\in GL_n(R)$). Then $\theta(T)\in E_n(R[T])$.
\et

%\rmk
%Using a result of Popescu \cite{po,sw}, the above result of Vorst can be extended to the case 
%when $R$ is a regular ring containing a field.

\subsection{Homotopy orbits of unimodular rows}

For a ring $R$, consider the set 
$$Um_{n+1}(R):=\{(a_1,\cdots,a_{n+1})\in R^{n+1}\,|\, \sum_{i=1}^{n+1} a_i b_i=1 \text{ for some } b_1,\cdots,b_{n+1}\in R\}$$
of \emph{unimodular rows} of length $n+1$ in $R$.
 Two unimodular rows $(a_1,\cdots,a_{n+1})$ and 
$(a'_1,\cdots,a'_{n+1})$ are homotopic if there is $(f_1(T),\cdots,f_{n+1}(T))\in Um_{n+1}(R[T])$ such that
$f_i(0)=a_i$ and $f_i(1)=a'_i$ for $i=1,\cdots,n+1$. The set of equivalence classes with respect to the equivalence relation generated by homotopies will be denoted by $\pi_0(Um_{n+1}(R))$.

We shall need  the following theorem from \cite{dtz2} very soon. See also \cite[Theorem 2.1]{f1} for a more general version.

\bt\label{uni}\cite[2.3]{dtz2}
Let $R$ be a  regular ring which is essentially of finite type over a field $k$. Then, for any $n\geq 2$ 
there is a bijection  $\eta_R:\pi_0(Um_{n+1}(R))\stackrel{\sim}\lra Um_{n+1}(R)/E_{n+1}(R)$.
\et

\notation
Let $v=(a_1,\cdots,a_{d+1})\in Um_{d+1}(R)$. The  orbit of $v$ in   $ Um_{d+1}(R)/E_{d+1}(R)$ will be written as $[v]=[a_1,\cdots,a_{d+1}]$. 

%\bc\label{eqvuni}
%Let $R$ be a  regular ring which is essentially of finite type over a field $k$.  
%Then, for $n\geq 2$,
%the relation induced by homotopy on $Um_{n+1}(R)$ is an equivalence relation.
%\ec

%\proof
%Let $u,v\in Um_{n+1}(R)$ be in the same homotopy orbit. Then there is a finite sequence 
%$u=w_0, w_1,\cdots,w_{r-1}, w_{r}=v$ in $Um_{n+1}(R)$ such that $w_i$ is homotopic to $w_{i+1}$
%for $i=0,\cdots,r-1$. It follows from the above theorem that there is $\sigma_i\in E_{n+1}(R)$
%such that $w_i\sigma_i=w_{i+1}$, for $i=0,\cdots,r-1$.   In view of Remark \ref{elem},
%the product  $\prod_{i=0}^{r-1} \sigma_i\in E_{n+1}(R)$
%then yields the desired homotopy between $u$ and $v$. 
%\qed

 \subsection{The pointed set $Q_{2n}(R)$ and its homotopy orbits:}
Let $R$ be any commutative Noetherian ring. Let $n\geq 2$ and we recall 
the following set, which appeared in \cite{f}:  
$$Q_{2n}(R)=\{(x_1,\cdots x_n,y_1,\cdots,y_n,z)\in R^{2n+1}\,|\,\sum_{i=1}^{n} x_i y_i=z-z^2\}.$$
 By definition,  elements $(x_1,\cdots x_n,y_1,\cdots,y_n,z)$ and $(x'_1,\cdots x'_n,y'_1,\cdots,y'_n,z')$ of $Q_{2n}(R)$ are homotopic if there
is $(f_1,\cdots,f_n,g_1,\cdots,g_n,h)$ in $Q_{2n}(R[T])$ such that $f_i(0)=x_i$, $g_i(0)=y_i$,
$h(0)=z$, and $f_i(1)=x'_i$, $g_i(1)=y'_i$,
$h(1)=z'$. Consider the equivalence relation generated by homotopies on $Q_{2n}(R)$. The set of
equivalence classes will be denoted by $\pi_0(Q_{2n}(R))$.

\section{Generalities III: The  maps}\label{maps}
Let $R$ be a smooth affine domain of dimension $d\geq 2$ over an infinite perfect field $k$.
The purpose of this section is to define a set-theoretic map $\delta_R:E^d(R)\to Um_{d+1}(R)/E_{d+1}(R)$. This involves several steps.

\subsection{[The map $\theta_R:E^d(R)\lra\pi_{0}(Q_{2d}(R))$].}
We first recall the definition of a set-theoretic map from the Euler class group $E^d(R)$ to $\pi_0(Q_{2d}(R))$ from \cite{dtz1}.
By \cite[Remark 4.14]{brs1} we know that an arbitrary element of $E^d(R)$ can be represented by a single Euler cycle 
$(J,\omega_J)$, where $J$ is a reduced ideal of height $d$. Now $\omega_J:(R/J)^d\sur J/J^2$ is given by $J=\langle a_1,\cdots,a_d\rangle+J^2$, for
some $a_1,\cdots,a_d\in J$.
Applying the  Nakayama Lemma one obtains $s\in J^2$
such that  $J=\langle a_1,\cdots,a_d,s\rangle$ 
 with $s-s^2=a_1b_1+\cdots +a_d b_d$  for some $b_1,\cdots,b_d\in R$ (see \cite{mo1} for a proof).
We  associate to $(J,\omega_J)$ the homotopy class  $[(a_1,\cdots,a_d,b_1,\cdots,b_d,s)]$
in $\pi_0(Q_{2d}(R))$. 

In \cite[Proposition 4.2]{dtz1} we proved the following result. The reader may also 
consult \cite{af,mm} for a similar result proved using different methods than \cite{dtz1}.

\bp\label{map1}
Let $R$ be a regular domain of dimension $d\geq 2$ which is essentially of finite type over  an infinite perfect 
field $k$. 
The association $(J,\omega_J)\mapsto [(a_1,\cdots,a_d,b_1,\cdots,b_d,s)]$ is well defined and gives rise to a set-theoretic map $\theta_d: E^d(R)\to \pi_0(Q_{2d}(R))$. The map $\theta_d$ takes the
trivial Euler cycle to the homotopy orbit of the base point $(0,\cdots,0)$ of $Q_{2d}(R)$.
\ep

\subsection{[The map $\zeta_R:\pi_{0}(Q_{2d}(R))\lra Um_{d+1}(R)/E_{d+1}(R)]$.} The map we are about to define will again be a set-theoretic map. For a homotopy orbit $[(x_1,\cdots x_d,y_1,\cdots,y_d,z)]\in \pi_{0}(Q_{2d}(R)$, we assign 
$$\zeta_R([(x_1,\cdots x_d,y_1,\cdots,y_d,z)]):=[x_1,\cdots x_d,1-2z]$$

\bp
$\zeta_R:\pi_{0}(Q_{2d}(R))\lra Um_{d+1}(R)/E_{d+1}(R)$ is  well-defined.
\ep
\proof
We first note that 
$(x_1,\cdots,x_d,1-2z)\in Um_{d+1}(R)$.

By Thorem \ref{uni}, $\pi_0(Um_{d+1}(R))=Um_{d+1}(R)/E_{d+1}(R)$.
Although the homotopy is not an equivalence relation on $Q_{2d}(R)$, to check that $\zeta_R$ is well-defined, it is enough to show that if
$(x'_1,\cdots x'_d,y'_1,\cdots,y'_d,z')\in Q_{2d}(R)$ is homotopic to $(x_1,\cdots x_d,y_1,\cdots,y_d,z)$, then the unimodular rows $(x_1,\cdots x_d,1-2z)$ and $(x'_1,\cdots x'_d,1-2z')$
are homotopic. Let $(f_1,\cdots,f_d,g_1,\cdots,g_d,h)\in Q_{2d}(R[T])$ be such that  
$f_i(0)=x_i$, $g_i(0)=y_i$, $h(0)=z$, and $f_i(1)=x'_i$, $g_i(1)=y'_i$, $h(1)=z'$ ($1\leq i\leq d$).
Clearly, $(f_1,\cdots,f_d,1-2h)\in Um_{d+1}(R[T])$ gives the desired homotopy between 
the unimodular rows $(x_1,\cdots x_d,1-2z)$ and $(x'_1,\cdots x'_d,1-2z')$. \qed

\subsection{[The map $\delta_R:E^d(R)\lra Um_{d+1}(R)/E_{d+1}(R)$].}
Finally, the map $\delta_R$ is simply defined to be the composite:
$$E^d(R)\stackrel{\theta_R}{\lra} \pi_0(Q_{2d}(R))\stackrel{\zeta_R}{\lra}Um_{d+1}(R)/E_{d+1}(R)$$

Let us summarize the description of $\delta_R$. 
Let $(J,\omega_J)\in E^d(R)$, where $J$ is a reduced ideal of height $d$. Now $\omega_J:(R/J)^d\sur J/J^2$ is given by $J=\langle\,a_1,\cdots,a_d\rangle+J^2$, for
some $a_1,\cdots,a_d\in J$.
Applying the  Nakayama Lemma one obtains $s\in J^2$
such that  $J=\langle\,a_1,\cdots,a_d,s\rangle$ 
 with $s-s^2=a_1b_1+\cdots +a_d b_d$  for some $b_1,\cdots,b_d\in R$.   
$\delta_R$ 
 takes $(J,\omega_J)$  to the orbit 
$[a_1,\cdots a_d,1-2s]\in Um_{d+1}(R)/E_{d+1}(R)$. 

\smallskip

A series of remarks are in order.

\smallskip

\rmk
If the characteristic of $k$ is $2$, then clearly $\delta_R$ turns out to be the trivial map. 

\smallskip

\rmk\label{original}
In our orginial definition in \cite{dtz2}, we defined $\delta_R((J,\omega_J))=[2a_1,\cdots 2a_d,1-2s]$. Note that, if $\sqrt{2}\in R$, then the two definitions coincide (apply \cite[Lemma 3.5 (ii)]{vk2}). In particular, 
they do so when $k=\mathbb R$ and the results in \cite{dtz2} go through with the above definition of $\delta_R$.

\smallskip

\rmk
Note that $(1-2s)^2\equiv 1$ modulo the ideal
$\langle\,a_1,\cdots a_d\rangle$ and therefore, the image of $\delta_R$ is hitting the orbits of some special type of unimodular rows. 
Conversely,  
let an orbit $[v]=[x_1,\cdots,x_d,z]\in Um_{d+1}(R)/E_{d+1}(R)$
be such that the ideal $\langle\,x_1,\cdots,x_d\rangle$ is reduced of height $d$, and $z^2\equiv 1$ modulo  
$\langle\,x_1,\cdots,x_d\rangle$. If $\frac{1}{2}\in R$,  then $[v]$ is in the image of $\delta_R$.

\section{Special case: Mennicke-like group structure}\label{menn}
We will say that the group structure on $Um_{d+1}(R)/E_{d+1}(R)$ is \emph{Mennicke-like}\footnote{In literature it has been described as \emph{nice group structure}. Ravi Rao suggested us to use the term  \emph{Mennicke-like}.} if for two orbits $[a_1,\cdots,a_d,x],[a_1,\cdots,a_d,y]\in Um_{d+1}(R)/E_{d+1}(R)$ we have the \emph{coordinate-wise product}:
$$[a_1,\cdots,a_d,x][a_1,\cdots,a_d,y]=[a_1,\cdots,a_d,xy].$$

Throughout this section, let $R$ be a smooth affine domain of dimension $d\geq 2$ over an infinite perfect field $k$.

\bl\label{torsion}
Let the group structure on $Um_{d+1}(R)/E_{d+1}(R)$ be Mennicke-like. Let $(J,\omega_J)\in E^{d}(R)$ be any element. Then $\delta_R((J,\omega_J))$ is $2$-torsion.
\el

\proof
If  $Char(k)=2$, then $\delta_R$ is trivial and we are done. Therefore, we assume that $Char(k)\neq 2$.
Let $\omega_J$ be induced by $J=\langle\,a_1,\cdots,a_d\rangle+J^2$. Then, there exists $s\in J^2$ such that 
$J=\langle\,a_1,\cdots,a_d,s\rangle$ with $s-s^2\in \langle\,a_1,\cdots,a_d\rangle$.
By definition, $\delta_R((J,\omega_J))=[a_1,\cdots,a_d,1-2s]$. As the 
group law is Mennicke-like, 
$$[a_1,\cdots,a_d,1-2s]^2=[a_1,\cdots,a_d,(1-2s)^2]=[a_1,\cdots,a_d,1].\qed$$

\bt\label{mennicke}
Let the group structure on $Um_{d+1}(R)/E_{d+1}(R)$ be Mennicke-like. Then 
$\delta_R$ is a morphism of groups.
\et

\proof
As in the above lemma, we may assume that $Char(k)\neq 2$.
Let $(J,\omega_J), (K,\omega_K)\in E^d(R)$ be such that $J+K=R$, where $J,K$ are both reduced ideals of height $d$. 
Then $(J,\omega_J)+ (K,\omega_K)=(J\cap K,\omega_{J\cap K})$, where $\omega_{J\cap K}$ is induced by $\omega_J$ and $\omega_K$. 
To prove the theorem, it is enough to show that
$$\delta_R((J,\omega_J))\,\ast\,\delta_R((K,\omega_K))=\delta_R ((J\cap K,\omega_{J\cap K})),$$
where $\ast$ denotes the product in $Um_{d+1}(R)/E_{d+1}(R)$

Let $\omega_{J\cap K}$ be induced by 
$J\cap K=\langle\,a_1,\cdots,a_d\rangle+(J\cap K)^2$. Then $J=\langle\,a_1,\cdots,a_d\rangle+J^2$ and $K=\langle\,a_1,\cdots,a_d\rangle+K^2$. Let $J=\langle\,a_1,\cdots,a_d,s\rangle$ with $s-s^2\in \langle\,a_1,\cdots,a_d\rangle$ and $K=\langle\,a_1,\cdots,a_d,t\rangle$ with $t-t^2\in \langle\,a_1,\cdots,a_d\rangle$,
as usual. Then it follows that $J\cap K=\langle\,a_1,\cdots,a_d,st\rangle$ and $st-s^2t^2\in \langle\,a_1,\cdots,a_d\rangle$. 

By the definition of the map  $\delta_R$, we have:
\begin{enumerate}
\item
 $\delta_R ((J,\omega_J))=[a_1,\cdots,a_d,1-2s]$, 
\item
$\delta_R ((K,\omega_K))=[a_1,\cdots,a_d,1-2t]$, 
\item
$\delta_R((J\cap K,\omega_{J\cap K}))=[a_1,\cdots,a_d,1-2st]$. 
\end{enumerate}
As the group law in $Um_{d+1}(R)/E_{d+1}(R)$ is Mennicke-like, we have
$$[a_1,\cdots,a_d,1-2s][a_1,\cdots,a_d,1-2t]=[a_1,\cdots,a_d,1-2s-2t+4st].$$
Let us try to locate a pre-image of the element on the right hand side of the above equation. To this end, 
we consider the following ideal 
$$L=\langle\,a_1,\cdots,a_d,s+t-2st\rangle=\langle\,a_1,\cdots,a_d,s^2+t^2-2st\rangle=\langle\,a_1,\cdots,a_d,(s-t)^2\rangle$$
in $R$ and note that  $L+J\cap K=R$ (as $s-t$ is a unit modulo $J\cap K$). Let `bar' denote modulo $\langle a_1,\cdots,a_d\rangle$. Then, 
$$\ol{L} \cap \ol{J\cap K}=\langle\,\ol{s}\ol{t}\rangle\langle\,\ol{s}+\ol{t}-\ol{2st}\rangle=\langle\,\ol{s^2t}+\ol{st^2}-2\ol{s}\ol{t}\rangle=
\langle\,\ol{s}\ol{t}+\ol{s}\ol{t}-2\ol{s}\ol{t}\rangle=\langle\,\ol{0}\rangle,$$ 
and we have $L\cap (J\cap K)=\langle\,a_1,\cdots,a_d\rangle$. Therefore, $(L,\omega_L)+(J\cap K,\omega_{J\cap K})=0$, where
$\omega_L$ is induced by the images of $a_1,\cdots,a_d$ in $L/L^2$. It is easy to see that $\delta_R((L,\omega_L))=
[ a_1,\cdots,a_d,1-2s-2t+4st]$. Finally, we conclude (using (\ref{torsion})) that
$$\delta_R ((J,\omega_J))\,\ast\,\delta_R((K,\omega_K))=\delta_R((L,\omega_L))=\delta_R((J\cap K,\omega_{J\cap K}))^{-1}=\delta_R((J\cap K,\omega_{J\cap K})).\qed$$

\section{The general case}\label{general}
In this section treat the general case. Our line of arguments (Theorem \ref{main} aided by Proposition \ref{prep}) may be termed as \emph{``Mennicke-Newman for ideals"}. For the Mennicke-Newman Lemma for elementary orbits of unimodular rows, see \cite[Lemma 3.2]{vk3}.

\bl
Let $I_1,I_2$ be two comaximal ideals in a ring $R$ such that $I_1\neq I_1^2$ and $I_2\neq I_2^2$. Then we can find $x\in I_1\smallsetminus I_1^2$ and $y\in I_2\smallsetminus I_2^2$ such that $x+y=1$.
\el

\proof
As $I_1^2+I_2^2=R$, we can find $a\in I_1^2$, $b\in I_2^2$ such that $a+b=1$. 

Claim: $I_1\cap I_2\not\subseteq I_1^2$. To see this note that $I_1^2+I_2=R$, and we have
$$I_1=I_1\cap R=I_1\cap (I_1^2+I_2)=I_1^2+I_1\cap I_2.$$
If $I_1\cap I_2\subseteq I_1^2$, then $I_1=I_1^2$, contrary to the hypothesis. Similarly,
$I_1\cap I_2\not\subseteq I_2^2$. 

Therefore, we can choose $\alpha\in I_1\cap I_2\smallsetminus 
(I_1^2\cup I_2^2)$.
Take $x=a-\alpha$ and $y=b+\alpha$ to conclude.
\qed

\bp\label{prep}
Let $R$ be a ring of dimension $d\geq 2$.
Let $J=\mm_1\cap \cdots \cap\mm_r$ and $K=\mm_{r+1}\cap \cdots \cap \mm_s$ be two ideals, each of height $d$, where $\mm_i$ are all disinct maximal ideals for $i=1,\cdots ,s$. Then,
there exist $x\in J$ and $y\in K$ such that:
\begin{enumerate}
\item
$x+y=1$,
\item
$x\not\in \mm_1^2\cup \cdots \mm_s^2$ and $y\not\in \mm_1^2\cup \cdots \mm_s^2$. 

\end{enumerate}
\ep

\proof
As $J^2+K^2=R$, we can find $a\in J^2$ and $b\in K^2$ such that $a+b=1$. We claim that 
there exists $c\in J\cap K$ such that $c\not\in \mm_1^2\cup \cdots \mm_s^2$. If we can prove the claim, we will take $x=a-c$ and $y=a+c$ to prove the proposition.

\smallskip

\noindent
\emph{Proof of the claim.}  We have $\mm_1^2+\mm_2^2\cdots \mm_s^2=R$. Choose $f\in \mm_1^2$
and $g\in  \mm_2^2\cdots \mm_s^2$ so that $f+g=1$. 

Observe that $\mm_1\cap (\mm_2^2\cdots \mm_s^2)\not\subseteq \mm_1^2$ (to see this, use the above lemma to obtain $z\in \mm_1\smallsetminus \mm_1^2$ and $w\in \mm_2^2\cdots \mm_s^2$ so that
$z+w=1$. Assume, if possible, that 
$\mm_1\cap (\mm_2^2\cdots \mm_s^2)\subseteq \mm_1^2$. As $z=z^2+wz$ and $wz\in  \mm_1\cap (\mm_2^2\cdots \mm_s^2)$ it would follow that $z\in \mm_1^2$. Contradiction.)

Choose $\alpha \in \mm_1\cap (\mm_2^2\cdots \mm_s^2)\smallsetminus \mm_1^2$ and take 
$c_1=f-\alpha$, $c_1'=g+\alpha$. Then, we have:  (1) $c_1+c_1'=1$, (2) $c_1\in \mm_1\smallsetminus \mm_1^2$, (3) $c_1\equiv 1$ modulo $\mm_i^2$ for all $i\neq 1$. 

Following a similar method, for each $i=1,\cdots ,s$,  choose $c_i\in \mm_i\smallsetminus \mm_i^2$ so that
$c_i\equiv 1$ modulo $\mm_1^2\cdots \mm_{i-1}^2\mm_{i+1}^2\cdots \mm_s^2$. Take 
$c=\prod_{i=1}^{s}c_i$. Then $c\in \mm_1\cdots \mm_s$ and it is easy to check that $c\not\in \mm_i^2$ for any 
$i$.   This completes the proof of the claim.
\qed

\bt\label{main}
Let $R$ be a smooth affine domain of dimension $d\geq 2$ over an infinite perfect field $k$. Then 
$\delta_R:E^d(R)\lra Um_{d+1}(R)/E_{d+1}(R)$ is a morphism of groups.
\et

\proof
If  $Char(k)=2$, then $\delta_R$ is trivial and we are done. Therefore, we assume that $Char(k)\neq 2$.

Let $(J,\omega_J), (K,\omega_K)\in E^d(R)$ be such that $J+K=R$, where $J,K$ are both reduced ideals of height $d$. 
Then $(J,\omega_J)+ (K,\omega_K)=(J\cap K,\omega_{J\cap K})$, where $\omega_{J\cap K}$ is induced by $\omega_J$ and $\omega_K$. 
To prove the theorem, it is enough to show that
$$\delta_R((J,\omega_J))\,\ast\,\delta_R((K,\omega_K))=\delta_R ((J\cap K,\omega_{J\cap K})),$$
where $\ast$ denotes the product in $Um_{d+1}(R)/E_{d+1}(R)$.

Let $J=\mathfrak{m}_1\cap \cdots \cap \mathfrak{m}_r$ and $K=\mathfrak{m}_{r+1}\cap \cdots \cap \mathfrak{m}_s$. Applying the above proposition, choose $x\in J$ and $y\in K$ such that 
 $x+y=1$
and $x\not\in \mm_1^2\cup \cdots \mm_s^2$ and $y\not\in \mm_1^2\cup \cdots \mm_s^2$.  Then $xy\in (J\cap K)\smallsetminus (J\cap K)^2$. As $x+y=1$, it is easy to check that for each $i$, the image of $xy$ in $\mathfrak{m}_i/\mathfrak{m}_i^2$ is not trivial. Therefore, $xy$ is a part of a basis of $\mathfrak{m}_i/\mathfrak{m}_i^2$, for each $i$,
$1\leq i\leq s$. Consequently, $xy$ is a part of generators of $(J\cap K)/(J\cap K)^2$.
 Similarly,  $x$ is a part of generators of $J/J^2$ and  $y$ is a part of generators of $K/K^2$. 

Let 
$J\cap K=\langle\,xy,a_1,\cdots,a_{d-1}\rangle+(J\cap K)^2$ for some $a_1,\cdots,a_{d-1}\in J\cap K$.
Let $\omega'_{J\cap K}:(R/(J\cap K))^d\sur (J\cap K)/(J\cap K)^2$ denote the corresponding surjection.  
 By \cite[2.2 and 5.0]{brs3} there is a unit $u$ modulo $J\cap K$ such that 
 $(J\cap K,\omega_{J\cap K})=(J\cap K,u\omega'_{J\cap K})$ in $E^d(R)$. 
Therefore, $(J\cap K,\omega_{J\cap K})$ is given by 
$J\cap K=\langle\,xy,ua_1,a_2,\cdots,a_{d-1}\rangle+(J\cap K)^2$. Similarly, 
$J=\langle\,x,ua_1,a_2,\cdots,a_{d-1}\rangle+J^2$ gives $(J,\omega_J)$ and 
$K=\langle\,x,ua_1,a_2,\cdots,a_{d-1}\rangle+K^2$ gives $(K,\omega_K)$. 

We can choose $s\in J\cap K$ such that $s-s^2\in \langle\,xy,ua_1,a_2,\cdots,a_{d-1}\rangle$ and 
$J\cap K=\langle\,xy,ua_1,a_2,\cdots,a_{d-1},s\rangle$. As 
$s-s^2\in \langle\,x,ua_1,a_2,\cdots,a_{d-1}\rangle$ and $s-s^2\in \langle\,y,ua_1,a_2,\cdots,a_{d-1}\rangle$ as well, it follows that $(J,\omega_J)$ corresponds to $J=\langle\,x,ua_1,a_2,\cdots,a_{d-1},s\rangle$ and $(K,\omega_K)$ corresponds to $K=\langle\,y,ua_1,a_2,\cdots,a_{d-1},s\rangle$ (to check this, use $x+y=1$). 

We then have,
\begin{enumerate} 
\item
$\delta_R((J,\omega_J))=[x,ua_1,a_2,\cdots,a_{d-1},1-2s],$
\item
$\delta_R((K,\omega_K))=[y,ua_1,a_2,\cdots,a_{d-1},1-2s],$
\item
$\delta_R ((J\cap K,\omega_{J\cap K}))=[xy,ua_1,a_2,\cdots,a_{d-1},1-2s].$
\end{enumerate} 

It follows that $\delta_R((J,\omega_J))\,\ast\,\delta_R((K,\omega_K))=\delta_R ((J\cap K,\omega_{J\cap K}))$, as $x+y=1$.
\qed

\section{A few remarks}
Let $R$ be a smooth affine domain of dimension $d\geq 2$ over an
infinite perfect field $k$. We now recall the definition of 
a group homomorphism $\phi_R:Um_{d+1}(R)/E_{d+1}(R)\to E^d(R)$. When $d$ is even,
$\phi_R$ has been defined in \cite{brs3}. The extension to general $d$ is available in \cite{dz,vk4}.
We urge the reader to look at \cite[Section 4]{dz} for the details.

\smallskip

\definition\label{phi}
Let $v=(a_1,\cdots,a_{d+1})\in Um_{d+1}(R)$.  Applying elementary transformations if necessary, we may assume that the height of the ideal $\langle a_1,\cdots,a_d\rangle$ is $d$. Write $J=\langle a_1,\cdots,a_d\rangle$ and let
$\omega_J:R^d\sur J$ be the surjection induced by $a_1,\cdots,a_d$. As $a_{d+1}$ is a unit modulo $J$,
we have $J=\langle a_1,\cdots,a_da_{d+1}\rangle+J^2$ and the corresponding element in $E^d(R)$ is $(J,\ol{a_{d+1}}\omega_J)$. Let $[v]$ denote the orbit of $v$ in $Um_{d+1}(R)/E_{d+1}(R)$.  Define $\phi_R([v])=(J,\ol{a_{d+1}}\omega_J)$. It is proved in \cite{dz,vk4} that $\phi_R$ is a morphism.

In passing, we record some observations on the composite maps  $\delta_R\phi_R$ and $\phi_R\delta_R$
(the latter played a crucial role in \cite{dtz2}).

\bt\label{comp1}
Let $R$ be a smooth affine domain of dimension $d\geq 2$ over an
infinite perfect field $k$.
For any  $(J,\omega_J)\in E^d(R)$, we have
$$\phi_R\delta_R((J,\omega_J))=(J,\omega_J)-(J,-\omega_J).$$
\et

\proof
Essentially the same proof as \cite[Theorem 2.10]{dtz2}. \qed

\bt\label{comp2}
Let $R$ be a smooth affine domain of dimension $d\geq 2$ over an
infinite perfect field $k$. Assume further that $\sqrt{2}\in R$.
Let $v=(a_1,\cdots,a_{d+1})\in  Um_{d+1}(R)$
 and let $v^{\ast}$ denote the ``antipodal" vector
$(-a_1,\cdots,a_{d+1})\in  Um_{d+1}(R)$. Then we have
$$\delta_R\phi_R([v])=[v][{v^{\ast}}]^{-1}.$$
\et

\proof
Let $v=(a_1,\cdots,a_{d+1})\in Um_{d+1}(R)$. Recall from (\ref{phi}) that $\phi_R$ takes $[v]$ to $(J,\omega_J)$ 
where $J=\langle a_1,\cdots,a_{d}\rangle$ and $\omega_J$ is induced by 
$J=\langle a_1,\cdots,a_{d-1},a_da_{d+1}\rangle +J^2.$
 Now as $v\in Um_{d+1}(R)$, there exist $b_1,\cdots,b_{d+1}\in R$ such that $a_1b_1+\cdots+a_{d+1}b_{d+1}=1$. Multiplying both sides by 
 $a_{d+1}b_{d+1}$ we get
 $$a_1b_1'+\cdots+a_{d-1}b_{d-1}'+a_da_{a_{d+1}}b_d'+(a_{d+1}b_{d+1})^2=a_{d+1}b_{d+1},$$
 implying that
$a_1b_1'+\cdots+a_{d-1}b_{d-1}'+a_da_{a_{d+1}}b_d'=a_{d+1}b_{d+1}-(a_{d+1}b_{d+1})^2.$
Note that $1-a_{d+1}b_{d+1}=a_1b_1+\cdots+a_{d}b_{d}\in J$. Therefore, 
$$\delta_R((J,\omega_J))=[a_1,\cdots,a_{d-1},a_da_{d+1},1-2(1-a_{d+1}b_{d+1})]$$
$$= [a_1,\cdots,a_{d-1},a_da_{d+1},2a_{d+1}b_{d+1}-1] =\delta_R\phi_R([v]).$$

Now consider $v^*=(a_1,\cdots,a_d,-a_{d+1})\in Um_{d+1}(R)$, the \emph{antipodal}  of $v$. Then by 
\cite[Lemma 3.5(iii)]{vk2}, 
$[a_1,\cdots,a_d,-a_{d+1}]^{-1}=[a_1,\cdots,a_d,b_{d+1}] \text{ in } Um_{d+1}(R)/E_{d+1}(R).$ 
 We now compute:
\begin{eqnarray*}
\delta_R\phi_R([v])=[a_1,\cdots,a_{d-1},a_da_{d+1},2a_{d+1}b_{d+1}-1]\\
=[a_1,\cdots,a_{d-1},2a_da_{d+1},2a_{d+1}b_{d+1}-1] \text{ (as 2 is a square) }\\
=[a_1,\cdots,a_{d+1}][a_1,\cdots,a_d,b_{d+1}] \text{ by \cite[3.5 (i)]{vk2} }\\
=[a_1,\cdots,a_{d+1}][a_1,\cdots,a_d,-a_{d+1}]^{-1} 
=[v][{v^{\ast}}]^{-1}.
\end{eqnarray*}
The proof is therefore complete.
\qed

\smallskip

Let $X=\Spec(R)$ be a smooth affine variety of dimension $d\geq 2$ over ${\mathbb R}$. Let $X({\mathbb R})$
denote the set of real points of $X$. Assume that $X({\mathbb R})\neq \emptyset$.
Then $X({\mathbb R})$ is a smooth real manifold. Let ${\mathbb R}(X)$ denote the ring obtained from $R$ by inverting all the functions which do not have any real zeros.
 We can apply  (\ref{mennicke}) to obtain the following result.

\bt
Let $X=\Spec(R)$ be a smooth affine variety of dimension $d\geq 2$ over ${\mathbb R}$ such that $X({\mathbb R})$ is orientable,
and the number of compact connected componenets of $X({\mathbb R})$ is at least one.
Then the group structure on $Um_{d+1}(R)/E_{d+1}(R)$ can never be Mennicke-like.
\et

\proof
In our earlier paper \cite{dtz2}, we proved the following assertions:
\begin{enumerate}
\item
$Um_{d+1}(R)/E_{d+1}(R)= Um_{d+1}({\mathbb R}(X))/E_{d+1}(({\mathbb R}(X))\bigoplus K$,
\item
$\delta_{{\mathbb R}(X)}:E^d({\mathbb R}(X))\lra Um_{d+1}({\mathbb R}(X))/E_{d+1}(({\mathbb R}(X))$ is an isomorphism,
\item
$Um_{d+1}({\mathbb R}(X))/E_{d+1}(({\mathbb R}(X))\iso \bigoplus_{t} \bz$, where $t$ is the number of compact connected components of $X({\mathbb R})$.
\end{enumerate}

 Now, assume that  the group structure on 
$Um_{d+1}({\mathbb R}(X))/E_{d+1}({\mathbb R}(X))$ is Mennicke-like. Then, by (\ref{mennicke}),
we shall find non-trivial orbits which are $2$-torsion. But as $t\geq 1$, the group $Um_{d+1}({\mathbb R}(X))/E_{d+1}({\mathbb R}(X))$ is non-trivial and is free abelian. Thus we arrive at a contradiction. As $Um_{d+1}({\mathbb R}(X))/E_{d+1}({\mathbb R}(X))$ is a subgroup of $Um_{d+1}(R)/E_{d+1}(R)$, the theorem follows.
\qed

\rmk
In \cite{dtz2} we also computed the universal Mennicke symbol $MS_{d+1}(R)$, where $R$ is as in the above theorem. It follows from there as well that the group structure on $Um_{d+1}(R)/E_{d+1}(R)$ can never be Mennicke-like. The arguments given above only avoids the computation of $MS_{d+1}(R)$.

\medskip

We now comment on the case when $X({\mathbb R})$ is non-orientable.

\bt
Let $X=\Spec(R)$ be a smooth affine variety of dimension $d\geq 2$ over ${\mathbb R}$ such that $X({\mathbb R})$ is 
non-orientable. Then $\delta_{{\mathbb R}(X)}:E^d({\mathbb R}(X))\lra Um_{d+1}({\mathbb R}(X))/E_{d+1}({\mathbb R}(X))$ is a surjective morphism. As a consequence, $Um_{d+1}({\mathbb R}(X))/E_{d+1}({\mathbb R}(X))$ is an $\bz/2\bz$-vector space of dimension $\leq t$, where $t$ is the number of compact connected componenets of $X({\mathbb R})$. 
\et

\proof
It has already been proved in \cite[Theorem 3.2]{dtz2} that $\delta_{{\mathbb R}(X)}$ is surjective.
In this article we proved that $\delta_{{\mathbb R}(X)}$ is a morphism. As $E^d({\mathbb R}(X))=\bigoplus_t \bz/2\bz$, the result follows. \qed

\end{document}